\documentclass[12pt]{article}

\usepackage{amsthm,amsmath,amssymb,latexsym}

\newtheorem{thm}{Theorem}[section]

\newtheorem{cor}[thm]{Corollary}

\theoremstyle{remark}
\newtheorem*{rmk}{Remark}

\numberwithin{equation}{section}
\renewcommand{\qed}{{\hfill\rule{4pt}{7pt}}\medskip}

\begin{document}
\begin{center}
{\Large\bf  New Finite Rogers-Ramanujan Identities}
\end{center}
\vskip 2mm
\centerline{Victor J. W. Guo$^1$, Fr\'ed\'eric Jouhet$^2$ and Jiang Zeng$^3$}

\begin{center}
{\footnotesize $^1$Department of Mathematics, East China Normal University, Shanghai 200062,
 People's Republic of China\\
{\tt jwguo@math.ecnu.edu.cn,\quad http://math.ecnu.edu.cn/\textasciitilde{jwguo}}\\[10pt]

$^2$Universit\'e de Lyon,
Universit\'e Lyon 1,
UMR 5208 du CNRS, Institut Camille Jordan,\\
F-69622, Villeurbanne Cedex, France\\
{\tt jouhet@math.univ-lyon1.fr,\quad http://math.univ-lyon1.fr/\textasciitilde{jouhet}}\\[10pt]

$^3$Universit\'e de Lyon,
Universit\'e Lyon 1,
UMR 5208 du CNRS, Institut Camille Jordan,\\
F-69622, Villeurbanne Cedex, France\\
{\tt zeng@math.univ-lyon1.fr,\quad http://math.univ-lyon1.fr/\textasciitilde{zeng}}
}
\end{center}

%\htmladdnormallink{http://math.univ-lyon1.fr/$\sim$zeng}{http://math.univ-lyon1.fr/~zeng}

\vskip 0.7cm {\small \noindent{\bf Abstract.}
We present  two general finite extensions for each of the two
Rogers-Ramanujan identities. Of these one can be derived directly
from Watson's transformation formula  by
specialization or through Bailey's method, the second
similar formula can be proved either by using the first formula and
the $q$-Gosper algorithm, or through the so-called Bailey lattice.
}


\section{Introduction}

The famous Rogers-Ramanujan identities (see \cite{Andrews76}) may be stated as follows:
\begin{align}
1+\sum_{k=1}^\infty\frac{q^{k^2}}{(1-q)(1-q^2)\cdots (1-q^k)}
&=\prod_{n=0}^\infty\frac{1}{(1-q^{5n+1})(1-q^{5n+4})},  \label{eq:rr-1}\\
1+\sum_{k=1}^\infty\frac{q^{k^2+k}}{(1-q)(1-q^2)\cdots (1-q^k)}
&=\prod_{n=0}^\infty\frac{1}{(1-q^{5n+2})(1-q^{5n+3})}. \label{eq:rr-2}
\end{align}
Throughout this paper we suppose that $q$ is a complex number such
that $0<|q|<1$. For any integer $n$ define the \emph{$q$-shifted factorial} $(a)_n$ by
$(a)_0=1$ and
$$
(a)_n=(a;q)_n:=
\begin{cases}
    %1, & \text{$n=0,$} \\
    (1-a)(1-aq)\cdots (1-aq^{n-1}), & \text{$n=1,2,\ldots,$} \\
    \left[(1-aq^{-1})(1-aq^{-2})\cdots (1-aq^{n})\right]^{-1}, & \text{$n=-1,-2,\ldots.$}
\end{cases}
$$
Note that  $\frac{1}{(q)_n}=0$ if $n<0$. For $m\geq 1$ we will also use
the compact notations:
\begin{align*}
 (a_1,\ldots,a_m)_n:=(a_1)_n\cdots (a_m)_n,
\qquad
(a_1,\ldots,a_m)_\infty:=\lim_{n\to\infty}(a_1,\ldots,a_m)_n.
\end{align*}

The following finite forms of the Rogers-Ramanujan
identities were proposed by Andrews~\cite{Andrews70}
(see also \cite{Andrews74,Bressoud81a})
as a shortcut to the Rogers-Ramanujan
identities:
\begin{align}
\sum_{k=0}^\infty\frac{q^{k^2}}{(q)_k (q)_{n-k}}
&=\sum_{k=-\infty}^\infty\frac{(-1)^k q^{(5k^2-k)/2}}{(q)_{n-k}(q)_{n+k}},\label{eq:andrews1}\\
\sum_{k=0}^\infty\frac{q^{k^2+k}}{(q)_k (q)_{n-k}}
&=\sum_{k=-\infty}^\infty\frac{(-1)^k q^{(5k^2-3k)/2}}{(q)_{n-k}(q)_{n+k}},\label{eq:andrews2}
\end{align}
where  $n$ is a nonnegative integer. These identities were first
proved by Watson~\cite{Watson} from his $q$-analogue of
Whipple's transformation formula and were
first stated in this form by Andrews. Later
Ekhad and Tre~\cite{ET} and Paule~\cite{Paule94} also used \eqref{eq:andrews1}
 as  the starting point of
the computer proofs of
the first Rogers-Ramanujan identity~\eqref{eq:rr-1}.

In the literature there are many  other
finite forms of the Rogers-Ramanujan identities due to
Schur~\cite{Schur}, Watson \cite{Watson}, Andrews
\cite{Andrews70}, Bressoud~\cite{Bressoud81b}, and
Warnaar~\cite{Warnaar}. We refer the reader to Sills \cite{Sills}
and Berkovich and Warnaar \cite{BW} for recent developments of this
subject.

The Bailey lemma~\cite{Bailey48} provides one of the most effective methods for proving
$q$-series identities of Rogers-Ramanujan type.
The main purpose of this paper is to study further
generalizations of Andrews' identities \eqref{eq:andrews1} and \eqref{eq:andrews2}
by using
Bailey's lemma, Bailey's chain and lattice, though we provide also alternative proofs.

First of all it is easy to derive the following general identities
 from Watson's
classical $q$-Whipple transformation
\cite[(III.18)]{GR}.

\begin{thm}\label{thm:lmnrs12}
If at least one of $a$, $b$ and $c$ is of the form
$q^n$, $n=1,2,\ldots,$ then the following identities hold
\begin{align}
&\hskip -3mm \sum_{k=-\infty}^\infty
\frac{(q/a,q/b,q/c,q/d,q/e)_k}{(a,b,c,d,e)_k}(abcdeq^{-3})^k
\nonumber\\
&\hspace{2cm}=\frac{(q,ab/q,bc/q,ac/q)_\infty}{(a,b,c,abc/q^2)_\infty}
\sum_{k=0}^\infty\frac{(q/a,q/b,q/c,de/q)_k}{(q,q^3/abc,d,e)_k}q^k,
\label{eq:abcde1}
\end{align}
\begin{align}
&\hskip -3mm \sum_{k=-\infty}^\infty
\frac{(q/a,q/b,q/c,q/d,q/e)_k}{(aq,bq,cq,dq,eq)_k}(abcdeq^{-1})^k
\nonumber\\
&\hspace{2cm}=\frac{(q,ab,bc,ac)_\infty}{(aq,bq,cq,abc/q)_\infty}
\sum_{k=0}^\infty\frac{(q/a,q/b,q/c,de)_k}{(q,q^2/abc,dq,eq)_k}q^k.
\label{eq:abcde2}
\end{align}
\end{thm}
Indeed,
by collecting the terms indexed by $k$ and $-k$ (resp. $-k-1$) together,
we see that \eqref{eq:abcde1} (resp. \eqref{eq:abcde2}) is in fact
 the limiting case $a\to 1$ (resp. $a\to q$)
of Watson's classical $q$-Whipple transformation
\cite[(III.18)]{GR} up to some obvious parameter replacements.

Note that, letting $d,e\to 0$ in \eqref{eq:abcde1} and
\eqref{eq:abcde2}, we obtain two simpler identities.
\begin{cor} We have
\begin{align}
\sum_{k=-\infty}^\infty \frac{(q/a,q/b,q/c)_k}{(a,b,c)_k}(abc)^k
q^{k^2-2k}
&=\frac{(q,ab/q,bc/q,ac/q)_\infty}{(a,b,c,abc/q^2)_\infty}
  \sum_{k=0}^\infty\frac{(q/a,q/b,q/c)_k}{(q,q^3/abc)_k}q^k, \label{eq:liu1}\\
\sum_{k=-\infty}^\infty
\frac{(q/a,q/b,q/c)_k}{(aq,bq,cq)_k}(abc)^k q^{k^2}
&=\frac{(q,ab,bc,ac)_\infty}{(aq,bq,cq,abc/q)_\infty}
\sum_{k=0}^\infty\frac{(q/a,q/b,q/c)_k}{(q,q^2/abc)_k}q^k,
\label{eq:liu2}
\end{align}
provided that at least one of $a$, $b$ and $c$ is of the form
$q^n$, $n=1,2,\ldots.$
\end{cor}

In a recent paper~\cite[(7.16) and (7.24)]{Liu}
Liu stated \eqref{eq:liu1} and
\eqref{eq:liu2} {\it incorrectly} as non-terminating series.
For example, if $bc=q$, the left-hand side of \eqref{eq:liu1} becomes
\begin{align*}
\sum_{k=-\infty}^\infty \frac{(q/a)_k}{(a)_k}a^k
q^{k^2-k}=\frac{(q)_\infty}{(a)_\infty},
\end{align*}
while the right-hand side of \eqref{eq:liu1} is equal to $0$
(since $ab/q=1$). Similarly, if $bc=1$, the left-hand side of
\eqref{eq:liu2} becomes
\begin{align*}
\sum_{k=-\infty}^\infty \frac{(q/a)_k}{(aq)_k}a^k
q^{k^2}=\frac{(q)_\infty}{(aq)_\infty},
\end{align*}
while the right-hand side of \eqref{eq:liu2} is equal to $0$. Thus
\eqref{eq:liu1} and \eqref{eq:liu2} do not hold for
non-terminating series. Therefore Theorems~13 and 14 in \cite{Liu} are invalid too.
A valid version of the latter theorems will be given in this paper within the framework of
Bailey chain.

Since the identities \eqref{eq:abcde1} and \eqref{eq:abcde2} are
between two terminating series, setting $a=q^{n+1}$, $b=q^{l+1}$,
$c=q^{m+1}$, $d=q^{u+1}$ and $e=q^{v+1}$, we  can rewrite them,
respectively, as follows:
\begin{align}
&\hskip -3mm
\sum_{k=0}^\infty\frac{q^{k^2} (q)_{l+m+n-k}(q)_{u+v+k}}
{(q)_k (q)_{l-k}(q)_{m-k}(q)_{n-k} (q)_{u+k}(q)_{v+k}} \nonumber\\
&=\sum_{k=-\infty}^\infty\frac{(-1)^k q^{(5k^2-k)/2}
(q)_{l+m}{(q)_{l+n}(q)_{m+n}} (q)_u(q)_v(q)_{u+v}}
{(q)_{l-k}(q)_{m-k}(q)_{n-k}(q)_{u-k}(q)_{v-k}(q)_{l+k}(q)_{m+k}(q)_{n+k}
(q)_{u+k}(q)_{v+k}},            \label{eq:lmnrs1}
\end{align}
\begin{align}
&\hskip -3mm
\sum_{k=0}^\infty\frac{q^{k^2+k} (q)_{l+m+n-k+1}(q)_{u+v+k+1}}
{(q)_k (q)_{l-k}(q)_{m-k}(q)_{n-k}(q)_{u+k+1}(q)_{v+k+1}} \nonumber\\
&=\sum_{k=-\infty}^\infty\frac{(-1)^k q^{(5k^2+3k)/2}
(q)_{l+m+1}(q)_{m+n+1}(q)_{l+n+1} (q)_u(q)_v (q)_{u+v+1}}
{(q)_{l-k}(q)_{m-k}(q)_{n-k}(q)_{u-k}(q)_{v-k}
(q)_{l+k+1}(q)_{m+k+1}(q)_{n+k+1}(q)_{u+k+1}(q)_{v+k+1}},
\label{eq:lmnrs2}
\end{align}
where  $l,m,n,u,v$ are nonnegative integers.

Starting from the identities \eqref{eq:abcde1} and \eqref{eq:abcde2},
we can derive the following  new identities, which were originally found through
computer experiments.
\begin{thm}\label{thm:abcde34}
If at least one of $a$, $b$ and $c$ is of the form
$q^n$, $n=1,2,\ldots,$ then the following identities hold
\begin{align}
&\hskip -3mm \sum_{k=-\infty}^\infty
\frac{(q/a,q/b,q/c,q/d,q/e)_k}{(a,b,c,d/q,e/q)_k}(abcdeq^{-3})^k
\nonumber\\
&=\frac{(q,ab/q,bc/q,ac/q)_\infty}{(a,b,c,abc/q^2)_\infty}
\sum_{k=0}^\infty\frac{(q/a,q/b,q/c,de/q^2)_k}{(q,q^3/abc,d,e)_k}q^k,
\label{eq:abcde3}
\end{align}
and
\begin{align}
&\hskip -3mm \sum_{k=-\infty}^\infty
\frac{(q/a,q/b,q/c,q/d,q/e)_k}{(a,b,c,d/q,e/q)_k}(abcdeq^{-4})^k
\nonumber\\
&=\frac{(q,ab/q,bc/q,ac/q)_\infty}{(a,b,c,abc/q^2)_\infty}
\sum_{k=0}^\infty\frac{(q/a,q/b,q/c,de/q^2)_k}{(q,q^3/abc,d,e)_k}q^{2k}.
\label{eq:abcde4}
\end{align}
\end{thm}

Setting $a=q^{n+1}$, $b=q^{l+1}$, $c=q^{m+1}$, $d=q^{u+1}$ and
$e=q^{v+1}$, where $l,m,n,u,v\in\mathbb{N}$, the above two identities
 can be written, respectively,  as follows:
\begin{align}
&\hskip -3mm
\sum_{k=0}^\infty\frac{q^{k^2} (q)_{l+m+n-k}(q)_{u+v+k-1}}
{(q)_k (q)_{l-k}(q)_{m-k}(q)_{n-k} (q)_{u+k}(q)_{v+k}} \nonumber\\
&=\sum_{k=-\infty}^\infty\frac{(-1)^k q^{(5k^2-k)/2}
(q)_{l+m}{(q)_{l+n}(q)_{m+n}} (q)_{u-1}(q)_{v-1}(q)_{u+v-1}}
{(q)_{l-k}(q)_{m-k}(q)_{n-k}(q)_{u-k}(q)_{v-k}(q)_{l+k}(q)_{m+k}(q)_{n+k}
(q)_{u+k-1}(q)_{v+k-1}},
\label{eq:lmnrs3}
\end{align}
and
\begin{align}
&\hskip -3mm
\sum_{k=0}^\infty\frac{q^{k^2+k} (q)_{l+m+n-k}(q)_{u+v+k-1}}
{(q)_k (q)_{l-k}(q)_{m-k}(q)_{n-k} (q)_{u+k}(q)_{v+k}} \nonumber\\
&=\sum_{k=-\infty}^\infty\frac{(-1)^k q^{(5k^2-3k)/2} (q)_{l+m}{(q)_{l+n}(q)_{m+n}}
(q)_{u-1}(q)_{v-1}(q)_{u+v-1}}
{(q)_{l-k}(q)_{m-k}(q)_{n-k}(q)_{u-k}(q)_{v-k}(q)_{l+k}(q)_{m+k}(q)_{n+k}
(q)_{u+k-1}(q)_{v+k-1}}.                           \label{eq:lmnrs4}
\end{align}

It is clear that letting $l,m,u,v\to\infty$
in \eqref{eq:lmnrs1} and \eqref{eq:lmnrs3}
(resp. \eqref{eq:lmnrs2} and
\eqref{eq:lmnrs4})
we obtain \eqref{eq:andrews1} (resp. \eqref{eq:andrews2}).

In Section~2 we will show connections between
 Theorems~\ref{thm:lmnrs12} and \ref{thm:abcde34} and the Bailey chain concept.
 In particular, we shall give a valid version
 (cf. Corollaries~2.2 and 2.3 below) of Liu's Theorems 13 and 14 by
iterating the Bailey lemma.
Since the left-hand sides of \eqref{eq:abcde3} and \eqref{eq:abcde4}
are \emph{not well-poised}, a modified version of Bailey's lemma, which is
called the Bailey lattice, is needed to prove \eqref{eq:abcde3}
while \eqref{eq:abcde4} is simply a combination of the previous ones.
We will give two proofs of \eqref{eq:abcde3} in Section~3 and
the proof of \eqref{eq:abcde4} in Section~\ref{sec:abcde4}. Finally we
derive  some interesting special cases  of Theorems~\ref{thm:lmnrs12}
and~\ref{thm:abcde34} in Section~\ref{sec:special}.

%%%%%%%%%%%%%%%%%%%%%%%%%%%%%%%%%%%%%%%%%%%%%%%%%%%%%%%%%%%%%%
\section{Bailey's lemma revisited}
%%%%%%%%%%%%%%%%%%%%%%%%%%%%%%%%%%%%%%%%%%%%%%%%%%%%%%%%%%%%%
The following theorem can  be found in~\cite[Chapter 12]{AAR}
(see also \cite{An86,Paule87}). Its original form is due to Bailey~\cite{Bailey48}.
The full power of the Bailey lemma  was discovered by Andrews
when he observed the iterative nature of the lemma, leading to the Bailey chain.
\begin{thm}[Bailey's lemma]
If $(\alpha_n, \beta_n)$ is a Bailey pair with parameter $x$, i.e., they are related by
\begin{equation}\label{pair}
\beta_n=\sum_{r\geq 0}{\alpha_r\over (q)_{n-r}(xq)_{n+r}},\qquad
\forall\ n\geq 0,
\end{equation}
then $(\alpha_n', \beta_n')$ is also a Bailey pair with parameter $x$,
where
\begin{align*}
\alpha_n'&=\frac{(\rho_1,\rho_2)_n(xq/\rho_1\rho_2)^n}
{(xq/\rho_1)_k(xq/\rho_2)_n}\alpha_n,\\
\beta_n'
&=\sum_{r\geq 0}\frac{(\rho_1,\rho_2)_r(xq/\rho_1\rho_2)_{n-r}(xq/\rho_1\rho_2)^r}
{(q)_{n-r}(xq/\rho_1)_n(xq/\rho_2)_n}\beta_r.
\end{align*}
\end{thm}
Thus we can iterate the above procedure to produce the so-called Bailey chain:
$$
(\alpha_n,\beta_n)\longrightarrow (\alpha_n',\beta_n')\longrightarrow (\alpha_n'',\beta_n'')\longrightarrow \cdots
$$
When $x=1$ we can symmetrize Bailey's lemma as follows.
\begin{cor}
If  two sequences $(\alpha_n)$ and $(\beta_n)$ are related by
\begin{equation*}%\label{pair'}
\beta_n=\sum_{r=-\infty}^\infty{\alpha_r\over (q)_{n-r}(q)_{n+r}},\qquad
\forall\ n\geq 0,
\end{equation*}
then
\begin{align}
&\sum_{n=-\infty}^\infty\frac{(\rho_1,\rho_2, q^{-N})_n}{(q/\rho_1, q/\rho_2, q^{N+1})_n}
\left(\frac{q^{1+N}}{\rho_1\rho_2}\right)^n(-1)^nq^{-{n\choose 2}}\alpha_n\nonumber\\
&\hspace{3cm}=\frac{(q, q/\rho_1\rho_2)_N}{(q/\rho_1, q/\rho_2)_N}
\sum_{n\geq 0}\frac{(\rho_1, \rho_2, q^{-N})_nq^n\beta_n}{(\rho_1\rho_2q^{-N})_n}.\label{x=1}
\end{align}
\end{cor}
\begin{proof}
Note that
the identity~\eqref{pair} for $(\alpha_n', \beta_n')$ can be written as
\begin{align}\label{equiv}
&\sum_{n\geq 0}\frac{(\rho_1,\rho_2, q^{-N})_n}{(xq/\rho_1, xq/\rho_2, xq^{N+1})_n}
\left(\frac{xq^{1+N}}{\rho_1\rho_2}\right)^n(-1)^nq^{-{n\choose 2}}\alpha_n\nonumber\\
&\hspace{2 cm}=\frac{(xq, xq/\rho_1\rho_2)_N}{(xq/\rho_1, xq/\rho_2)_N}
\sum_{n\geq 0}\frac{(\rho_1, \rho_2, q^{-N})_nq^n\beta_n}{(\rho_1\rho_2q^{-N}/x)_n}.
\end{align}
So, setting $x=1$ and
replacing  $\alpha_n$ by $\alpha_n+\alpha_{-n}$ for $n\geq 1$ in \eqref{equiv}
 yields \eqref{x=1}.
\end{proof}
\begin{rmk}
Set $A_n=q^n(q)_{2n}\beta_n$, $B_n=(-1)^nq^{-{n\choose 2}}\alpha_n$, $\rho_1=q/b$ and
$\rho_2=q/c$ in the above corollary
we obtain a valid version of  Theorem~13 in \cite{Liu} with $a=q^{1+N}$.
\end{rmk}

Let $\alpha_0=\beta_0=1$ and
$$
\alpha_n=(-1)^n(q^{n\choose 2}+q^{-n\choose 2}),\qquad \beta_n=0,\qquad n\geq 1.
$$
Then the $q$-binomial formula implies that $(\alpha_n, \beta_n)$ is a Bailey pair
with parameter $x=1$. Namely,
\begin{equation*}%\label{eq:pair1}
\sum_{r=-\infty}^\infty\frac{(-1)^rq^{r\choose 2}}{(q)_{n-r}(q)_{n+r}}=\delta_{n,0}.
\end{equation*}
Moving to the right in the Bailey chain we get
\begin{align*}
\alpha_n'&=\frac{(\rho_1,\rho_2)_n(q/\rho_1\rho_2)^n}{(q/\rho_1, q/\rho_2)_n}(-1)^n
\left(q^{n\choose 2}+q^{-n\choose 2}\right),\\
\beta_N'&=\frac{(q/\rho_1\rho_2)_N}{(q,q/\rho_1, q/\rho_2)_N}.
\end{align*}
The corresponding identity~\eqref{pair} reads
\begin{equation*}%\label{eq:pair1'}
\sum_{n=-\infty}^\infty\frac{(\rho_1, \rho_2)_n(q/\rho_1\rho_2)^n(-1)^nq^{n\choose 2}}
{(q)_{N-n}(q)_{N+n}(q/\rho_1, q/\rho_2)_n}=
\frac{(q/\rho_1\rho_2)_N}{(q,q/\rho_1, q/\rho_2)_N}.
\end{equation*}
Now the next Bailey pair is
\begin{align*}
\alpha_n''&=\frac{(\rho_1, \rho_2,\rho_3,\rho_4)_n}{(q/\rho_1, q/\rho_2, q/\rho_3, q/\rho_4)_n}
\left(\frac{q^2}{\rho_1\rho_2\rho_3\rho_4}\right)^n (-1)^n
\left(q^{n\choose 2}+q^{-n\choose 2}\right),\\
\beta_N''&=\sum_{n\geq 0}\frac{(q/\rho_3\rho_4)_{N-n}
(\rho_3, \rho_4, q/\rho_1\rho_2)_n}{(q)_{N-n}(q/\rho_3,q/\rho_4)_N
(q, q/\rho_1, q/\rho_2)_n}\left(\frac{q}{\rho_3\rho_4}\right)^n.
\end{align*}
The corresponding identity~\eqref{pair} is nothing else
than \eqref{eq:abcde1} with $a=q^{1+N}$ by the following substitution:
$$
\rho_1=q/b,\, \rho_2=q/c,\, \rho_3=q/d,\, \rho_4=q/e.
$$
Similarly, when $x=q$, we can symmetrize Bailey's lemma as follows.

\begin{cor}
If  two sequences $(\alpha_n)$ and $(\beta_n)$ are related by
\begin{equation*}%\label{pair''}
\beta_n=\sum_{r=-\infty}^\infty{\alpha_r\over (q)_{n-r}(q)_{n+r+1}},\qquad
\forall\ n\geq 0,
\end{equation*}
then
\begin{align*}
&\sum_{n=-\infty}^\infty\frac{(\rho_1,\rho_2, q^{-N})_n}{(q^2/\rho_1, q^2/\rho_2, q^{N+2})_n}
\left(\frac{q^{2+N}}{\rho_1\rho_2}\right)^n(-1)^nq^{-{n\choose 2}}\frac{\alpha_n}{1-q}\nonumber\\
&\hspace{3cm }=\frac{(q^2, q^2/\rho_1\rho_2)_N}{(q^2/\rho_1, q^2/\rho_2)_N}
\sum_{n\geq 0}\frac{(\rho_1, \rho_2, q^{-N})_nq^n\beta_n}{(\rho_1\rho_2q^{-N-1})_n}.
\end{align*}
\end{cor}
\begin{proof}
In Bailey's lemma set $x=q$ and replace $\alpha_n$ by $\frac{1}{1-q}(\alpha_n+\alpha_{-n-1})$.
\end{proof}
\begin{rmk}
Set $A_n=q^n(q)_{2n+1}\beta_n$, $B_n=(-1)^nq^{-{n\choose 2}}\alpha_n$, $\rho_1=q/b$ and
$\rho_2=q/c$ in the above corollary
we obtain a valid version of Theorem~14 in \cite{Liu} with $a=q^{1+N}$.
\end{rmk}

Similarly, let
$$
\alpha_n=(-1)^n(q^{n\choose 2}+q^{-n-1\choose 2}),\qquad \beta_n=\delta_{n,0}.
$$
Then the $q$-binomial formula implies that
 $(\alpha_n,\beta_n)$ is a Bailey pair with parameter $x=q$,
i.e.,
$$
\sum_{r=-\infty}^\infty{(-1)^rq^{r\choose 2}\over (q)_{n-r}(q)_{n+r+1}}=\delta_{n,0}.
$$
Moving to the right in the Bailey chain we find
\begin{align*}
\alpha_n'&=\frac{(\rho_1,\rho_2)_n(q^2/\rho_1\rho_2)^n}{(q^2/\rho_1, q^2/\rho_2)_n}
\frac{(-1)^n}{1-q}
\left(q^{n\choose 2}+q^{-n-1\choose 2}\right),\\
\beta_N'&=\frac{(q^2/\rho_1\rho_2)_N}{(q,q^2/\rho_1, q^2/\rho_2)_N}.
\end{align*}
The corresponding identity~\eqref{pair} reads
\begin{equation*}%\label{eq:pair1}
\sum_{n=-\infty}^\infty\frac{(\rho_1, \rho_2)_n(q^2/\rho_1\rho_2)^n(-1)^nq^{n\choose 2}}
{(q)_{N-n}(q)_{N+n+1}(q^2/\rho_1, q^2/\rho_2)_n}=
\frac{(q^2/\rho_1\rho_2)_N}{(q,q^2/\rho_1, q^2/\rho_2)_N}.
\end{equation*}
Finally the next Bailey pair is
\begin{align*}
\alpha_n''&=\frac{(\rho_1, \rho_2,\rho_3,\rho_4)_n}{(q^2/\rho_1, q^2/\rho_2, q^2/\rho_3, q^2/\rho_4)_n}
\left(\frac{q^2}{\rho_1\rho_2\rho_3\rho_4}\right)^n \frac{(-1)^n}{1-q}
\left(q^{n\choose 2}+q^{-n-1\choose 2}\right),\\
\beta_N''&=\sum_{n\geq 0}\frac{(q^2/\rho_3\rho_4)_{N-n}
(\rho_3, \rho_4, q^2/\rho_1\rho_2)_n}{(q)_{N-n}(q^2/\rho_3,q^2/\rho_4)_N
(q, q^2/\rho_1, q^2/\rho_2)_n}\left(\frac{q^2}{\rho_3\rho_4}\right)^n.
\end{align*}
The corresponding identity~\eqref{pair} is nothing else than \eqref{eq:abcde2}
with $a=q^{1+N}$ by the following substitution:
$$
\rho_1=q/b,\, \rho_2=q/c,\, \rho_3=q/d,\, \rho_4=q/e.
$$

%%%%%%%%%%%%%%%%%%%%%%%%%%%%%%%%%%%%%%%%%%%%%%%%%%%%%%%%%%%%%%
\section{Two proofs of \eqref{eq:abcde3}}

%%%%%%%%%%%%%%%%%%%%%%%%%%%%%%%%%%%%%%%%%%%%%%%%%%%%%%%%%%%%%%%%%%%%%%%%%%%%
\subsection{First Proof of \eqref{eq:abcde3}} \label{sec:proof}
We shall prove the equivalent form \eqref{eq:lmnrs3}.
Writing
$$
(1-q^{M+N+2})(q)_{M}(q)_{N} =(q)_{M}(q)_{N+1}
+q^{N+1}(q)_{M+1}(q)_{N},
$$
with $M=l+m+n-k$ and $N=u+v+k-1$ we  see that the left-hand side
of \eqref{eq:lmnrs3} multiplied by $(1-q^{M+N+2})$ is equal to
\begin{align*}
L=\sum_{k=0}^\infty\frac{q^{k^2} (q)_{l+m+n-k}(q)_{u+v+k}} {(q)_k
(q)_{l-k}(q)_{m-k}(q)_{n-k} (q)_{u+k}(q)_{v+k}}
+\sum_{k=0}^\infty\frac{q^{k^2+k+u+v}
(q)_{l+m+n-k+1}(q)_{u+v+k-1}} {(q)_k (q)_{l-k}(q)_{m-k}(q)_{n-k}
(q)_{u+k}(q)_{v+k}}.
\end{align*}
Note that the first sum is exactly the left-hand side of
\eqref{eq:lmnrs1} while the second sum corresponds to the
left-hand side  of \eqref{eq:lmnrs2} with $u\to u-1$ and $v\to
v-1$. Combining  any two terms indexed by $k$ and $-k$ (respectively,
$-k-1$) in the right-hand side of \eqref{eq:lmnrs1} (respectively,
\eqref{eq:lmnrs2}) then yields
\begin{align*}
L=-\frac{(q)_{l+m}(q)_{l+n}(q)_{m+n}(q)_{u+v}}{(q)_l^2(q)_m^2(q)_n^2(q)_u(q)_v}
+\sum_{k=0}^\infty f_k,
\end{align*}
where
\begin{align*}
f_k&=\frac{(-1)^k q^{(5k^2-k)/2}(1+q^k)
(q)_{l+m}(q)_{l+n}(q)_{m+n}
(q)_u(q)_v(q)_{u+v}}{(q)_{l-k}(q)_{m-k}(q)_{n-k}(q)_{u-k}(q)_{v-k}
(q)_{l+k}(q)_{m+k}(q)_{n+k}(q)_{u+k}(q)_{v+k}}      \nonumber\\
&{}+\frac{(-1)^k q^{(5k^2+3k)/2+u+v} (1-q^{2k+1})
(q)_{l+m+1}(q)_{m+n+1}(q)_{l+n+1} (q)_{u-1}(q)_{v-1}
(q)_{u+v-1}}{(q)_{l-k}(q)_{m-k}(q)_{n-k}(q)_{u-k-1}(q)_{v-k-1}
(q)_{l+k+1}(q)_{m+k+1}(q)_{n+k+1}(q)_{u+k}(q)_{v+k}}.
%\label{eq:twosum2}
\end{align*}
On the other hand, combining the terms indexed with $k$ and $-k$
in the right-hand side of \eqref{eq:lmnrs3} multiplied by
$(1-q^{M+N+2})$ yields
\begin{align*}
R=-\frac{(q)_{l+m}(q)_{l+n}(q)_{m+n}(q)_{u+v-1}(1-q^{l+m+n+u+v+1})}
{(q)_l^2(q)_m^2(q)_n^2(q)_{u-1}(q)_{v-1}} +\sum_{k=0}^\infty g_k,
\end{align*}
where
\begin{align*}
g_k&=\frac{(-1)^k q^{(5k^2-k)/2}(q)_{l+m}(q)_{l+n}(q)_{m+n}
(q)_{u-1}(q)_{v-1}(q)_{u+v-1}}{(q)_{l-k}(q)_{m-k}(q)_{n-k}(q)_{u-k}(q)_{v-k}
(q)_{l+k}(q)_{m+k}(q)_{n+k}(q)_{u+k-1}(q)_{v+k-1}} \nonumber\\
&\quad{}\times(1-q^{l+m+n+u+v+1})
\left(1+\frac{q^k(q)_{u-k}(q)_{v-k}}{(q)_{u+k}(q)_{v+k}}\right).  %\label{eq:twosum3}
\end{align*}
Now applying the $q$-Gosper's algorithm we see  that
\begin{align}
f_k-g_k=F(k+1)-F(k), \label{eq:fkgk}
\end{align}
where
\begin{align*}
F(k)&=\frac{(-1)^k
q^{(5k^2-3k)/2+u+v}(1-q^{l+m+n+k+1})(q)_{l+m}(q)_{l+n}(q)_{m+n}
(q)_{u-1}(q)_{v-1}(q)_{u+v-1}}{(q)_{l-k}(q)_{m-k}(q)_{n-k}(q)_{u-k}(q)_{v-k}
(q)_{l+k}(q)_{m+k}(q)_{n+k}(q)_{u+k-1}(q)_{v+k-1}}.
\end{align*}
Summing \eqref{eq:fkgk} over $k\geq 0$ yields $L=R$. \qed

\begin{rmk}
Since
$$
\frac{q^{(5k^2-k)/2}}{(q)_{v+k-1}}
=\frac{q^{(5k^2-k)/2}}{(q)_{v+k}}-\frac{q^{(5k^2+k)/2}}{(q)_{v+k}}q^v,
$$
by \eqref{eq:lmnrs1} we have
\begin{align*}
&\hskip -3mm
\sum_{k=0}^\infty\frac{q^{k^2} (q)_{l+m+n-k}(q)_{u+v+k}}
{(q)_k (q)_{l-k}(q)_{m-k}(q)_{n-k} (q)_{u+k}(q)_{v+k}} \\
&=\sum_{k=-\infty}^\infty\frac{(-1)^k q^{(5k^2-k)/2} (q)_{l+m}{(q)_{l+n}(q)_{m+n}}
(q)_{u}(q)_{v-1}(q)_{u+v}}
{(q)_{l-k}(q)_{m-k}(q)_{n-k}(q)_{u-k}(q)_{v-k}(q)_{l+k}(q)_{m+k}(q)_{n+k}
(q)_{u+k}(q)_{v+k-1}}.
\end{align*}
\end{rmk}

\subsection{Proof of (\ref{eq:abcde3}) through the Bailey lattice}\label{lattice}
%%%%%%%%%%%%%%%%%%%%%%%%%%%%%%%%%%%%%%%%%%%%%%%%%%%%%%%%%%%%%%%%%%%%%%%%%%%%%%%%%%%%
Now it is possible to prove (\ref{eq:abcde3}) in the context of Bailey pairs,
but one needs a modified version of Bailey's lemma, which allows to
change the parameter $x$ between two iterations of Bailey's lemma.
This was defined as the Bailey lattice by Agarwal,
Andrews and Bressoud in \cite{AAB} (see also \cite{SW, Warnaar2}),
where they proved with this tool the full Andrews-Gordon identity:
\begin{thm}[Bailey lattice]\label{thm:balat}
If $(\alpha_n, \beta_n)$ is a Bailey pair with parameter $x$
then $(\alpha_n', \beta_n')$ is a Bailey pair with parameter $xq^{-1}$,
where
$\alpha'_0=1$ and for $n\geq 1$
$$\alpha_n'=(1-x)(x/\rho_1\rho_2)^n\frac{(\rho_1,\rho_2)_n}{(x/\rho_1,x/\rho_2)_n}
\left(\frac{\alpha_n}{1-xq^{2n}}-xq^{2n-2}\frac{\alpha_{n-1}}{1-xq^{2n-2}}\right),$$
and for $n\geq 0$
$$\beta_n'
=\sum_{r\geq 0}\frac{(\rho_1,\rho_2)_r(x/\rho_1\rho_2)_{n-r}(x/\rho_1\rho_2)^r}
{(q)_{n-r}(x/\rho_1)_n(x/\rho_2)_n}\beta_r.$$
\end{thm}
Consider for $n\geq 0$
\begin{equation}\label{eq:Baileylattice}
\alpha_n=(-1)^nq^{n\choose 2}\frac{1-q^{2n+1}}{1-q},\qquad \beta_n=\delta_{n,0}.
\end{equation}
Then by the $q$-binomial formula, $(\alpha_n,\beta_n)$ is a Bailey pair with
parameter $x=q$, and we can move to the right in the Bailey chain to get
the following Bailey pair with parameter $x=q$:
\begin{align*}
\alpha_n'&=\frac{(\rho_1,\rho_2)_n(q^2/\rho_1\rho_2)^n}{(q^2/\rho_1, q^2/\rho_2)_n}
(-1)^nq^{n\choose 2}\frac{1-q^{2n+1}}{1-q},\\
\beta_N'&=\frac{(q^2/\rho_1\rho_2)_N}{(q,q^2/\rho_1, q^2/\rho_2)_N}.
\end{align*}
In the next step, apply Theorem~\ref{thm:balat} to assert that $(\alpha''_n,\beta''_n)$
is a Bailey pair with parameter $x=1$, where
$\alpha''_0=1$ and for $n\geq 1$
$$\alpha_n''=(1-q)(q/\rho_3\rho_4)^n\frac{(\rho_3,\rho_4)_n}{(q/\rho_3,q/\rho_4)_n}
\left(\frac{\alpha'_n}{1-q^{2n+1}}-q^{2n-1}\frac{\alpha'_{n-1}}{1-q^{2n-1}}\right),$$
and for $N\geq 0$
$$\beta_N''
=\sum_{n\geq 0}\frac{(\rho_3,\rho_4)_n(q/\rho_3\rho_4)_{N-n}(q/\rho_3\rho_4)^n}
{(q)_{N-n}(q/\rho_3)_N(q/\rho_4)_N}\beta'_n.$$
Setting $\rho_1=q^2/d$, $\rho_2=q^2/e$, $\rho_3=q/b$ and $\rho_4=q/c$ yields
\begin{eqnarray}
\beta_N''&=&\sum_{n\geq 0}\frac{(q/b,q/c)_n(bc/q)_{N-n}(bc/q)^n}
{(q)_{N-n}(b,c)_N}\frac{(de/q^2)_n}{(q,d,e)_n}\nonumber\\
&=&\frac{(bc/q)_N}{(q,b,c)_N}\sum_{n\geq 0}\frac{(q^{-N},q/b,q/c,de/q^2)_n}
{(q,d,e,q^{2-N}/bc)_n}q^n\label{beta1}.
\end{eqnarray}
On the other hand, writing $\displaystyle\beta''_N=\sum_{n\geq 0}{\alpha''_n\over (q)_{N-n}(q)_{N+n}}$, we get
\begin{eqnarray}
\beta_N''&=&\frac{1}{(q,q)_N}\left[1+\sum_{n\geq 1}\frac{(q^{-N},q/b,q/c)_n}
{(q^{1+N},b,c)_n}(bcq^{N-1})^n\frac{(q^2/d,q^2/e)_n}{(d,e)_n}(deq^{-2})^n\right]\nonumber\\
&&\hskip 1cm+\frac{1}{(q,q)_N}\sum_{n\geq 1}\frac{(q^{-N},q/b,q/c)_n}
{(q^{1+N},b,c)_n}(bcq^{N-1})^n\frac{(q^2/d,q^2/e)_{n-1}}{(d,e)_{n-1}}q(deq^{-1})^{n-1}\nonumber\\
&=&\frac{1}{(q,q)_N}\sum_{n\geq 0}\frac{(q^{-N},q/b,q/c,q^2/d,q^2/e)_n}
{(q^{1+N},b,c,d,e)_n}(bcdeq^{N-3})^n\nonumber\\
&&\hskip 2cm+\frac{1}{(q,q)_N}\sum_{n\geq 1}\frac{(q^{-N},q/b,q/c,q/d,q/e)_n}
{(q^{1+N},b,c,d/q,e/q)_n}(bcdeq^{N-2})^n\nonumber\\
&=&\frac{1}{(q,q)_N}\sum_{n\in\mathbb{Z}}\frac{(q^{-N},q/b,q/c,q/d,q/e)_n}
{(q^{1+N},b,c,d/q,e/q)_n}(bcdeq^{N-2})^n\label{beta2},
\end{eqnarray}
where the last equality follows by replacing $n$ by $-n$ in the first sum.
Finally, equating (\ref{beta1}) and (\ref{beta2}) yields (\ref{eq:abcde3}),
where $a$ is replaced by $q^{1+N}$.

%%%%%%%%%%%%%%%%%%%%%%%%%%%%%%%%%%%%%%%%%%%%%%%%%%%%%%%
\subsection{Remarks}
%%%%%%%%%%%%%%%%%%%%%%%%%%%%%%%%%%%%%%%%%%%%%%%%%%%%%%%%%%%
In  \cite[Lemma 4.3]{SW}  Schilling and Warnaar proved another transformation,
different but closely related to Theorem~\ref{thm:balat}. If we apply the
latter formula to the same pair \eqref{eq:Baileylattice}, then we can derive the following
extension of (1.4), different from the previous ones:
\begin{align}\label{eq:abcde6-1}
&\hskip -3mm \sum_{k=-\infty}^\infty
\frac{(q/a,q/b,q/c,q/d,q/e)_k}{(a,bq,cq,dq,eq)_k}(abcdeq^{-1})^k
\nonumber\\
&\hspace{2cm}=\frac{(q,ab,bc,ac)_\infty}{(a,bq,cq,abc/q)_\infty}
\sum_{k=0}^\infty\frac{(q/a,q/b,q/c,de)_k}{(q,q^2/abc,dq,eq)_k}q^k,
\end{align}
where at least one of $a,b,c$ is of the form $q^n$, $n=1,2,\ldots.$

However, Eq.~\eqref{eq:abcde6-1} can be easily deduced from \eqref{eq:abcde2}. This is
because the left-hand side of \eqref{eq:abcde6-1} can be written as
\begin{align*}
\sum_{k=-\infty}^\infty
\frac{(1-aq^k)(q/a,q/b,q/c,q/d,q/e)_k }{(1-a)(aq,bq,cq,dq,eq)_k }(abcdeq^{-1})^k,
\end{align*}
and we have
\begin{align}\label{eq:abcde=0}
\sum_{k=-\infty}^\infty
\frac{(q/a,q/b,q/c,q/d,q/e)_k }{(aq,bq,cq,dq,eq)_k }(abcde)^k=0
\end{align}
since
\begin{align}\label{eq:akinv}
\frac{(q/a)_{-k-1}}{(aq)_{-k-1}}a^{-k-1}
=\frac{(1/a)_{k+1}}{(a)_{k+1}}a^{k+1}
=-\frac{(q/a)_{k}}{(aq)_{k}}a^{k}.
\end{align}

An identity similar to \eqref{eq:abcde6-1} is as follows:
\begin{align}\label{eq:abcde6-3}
&\hskip -3mm \sum_{k=-\infty}^\infty
\frac{(q/a,q/b,q/c,q/d,q/e)_k}{(a,bq,cq,dq,eq)_k}(abcde)^k
\nonumber\\
&\hspace{2cm}=\frac{a(q,ab,bc,ac)_\infty}{q(a,bq,cq,abc/q)_\infty}
\sum_{k=0}^\infty\frac{(q/a,q/b,q/c,de)_k}{(q,q^2/abc,dq,eq)_k}q^k,
\end{align}
where at least one of $a,b,c$ is of the form $q^n$, $n=1,2,\ldots.$

Eq.\eqref{eq:abcde6-3} is also a consequence of \eqref{eq:abcde2}, because
the left-hand side of
\eqref{eq:abcde6-3} may be written as
\begin{align*}
\sum_{k=-\infty}^\infty
\frac{(1-aq^k)(q/a,q/b,q/c,q/d,q/e)_k}{(1-a)(aq,bq,cq,dq,eq)_k}(abcde)^k.
\end{align*}
But we have \eqref{eq:abcde=0} and, by \eqref{eq:akinv},
Eq.~\eqref{eq:abcde2} may be rewritten as
\begin{align*}
&\hskip -3mm \sum_{k=-\infty}^\infty
\frac{(q/a,q/b,q/c,q/d,q/e)_k}{(aq,bq,cq,dq,eq)_k}(abcdeq)^k
\nonumber\\
&\hspace{2cm}=-\frac{(q,ab,bc,ac)_\infty}{q(aq,bq,cq,abc/q)_\infty}
\sum_{k=0}^\infty\frac{(q/a,q/b,q/c,de)_k}{(q,q^2/abc,dq,eq)_k}q^k.
%\label{eq:abcde2}
\end{align*}

Note that \eqref{eq:abcde6-3} is an extension of
\begin{align*}
\sum_{k=0}^\infty\frac{q^{k^2+k}}{(q)_k (q)_{n-k}}
&=q^{-n}\sum_{k=-\infty}^\infty\frac{(-1)^k q^{(5k^2-5k)/2}}{(q)_{n-k}(q)_{n+k}}.
\end{align*}

Similarly to \eqref{eq:abcde6-1} and \eqref{eq:abcde6-3},
we can prove
\begin{align}\label{eq:abcde6-2}
&\hskip -3mm \sum_{k=-\infty}^\infty
\frac{(q/a,q/b,q/c,q/d,q/e)_k}{(aq,bq,cq,dq,e)_k}(abcdeq^{-1})^k
\nonumber\\
&\hspace{2cm}=\frac{(q,ab,bc,ac)_\infty}{(1-e)(aq,bq,cq,abc/q)_\infty}
\sum_{k=0}^\infty\frac{(q/a,q/b,q/c,de)_k}{(q,q^2/abc,dq,eq)_k}q^k, \\
\label{eq:abcde6-4}
&\hskip -3mm \sum_{k=-\infty}^\infty
\frac{(q/a,q/b,q/c,q/d,q/e)_k}{(aq,bq,cq,dq,e)_k}(abcde)^k
\nonumber\\
&\hspace{2cm}=\frac{e(q,ab,bc,ac)_\infty}{q(1-e)(aq,bq,cq,abc/q)_\infty}
\sum_{k=0}^\infty\frac{(q/a,q/b,q/c,de)_k}{(q,q^2/abc,dq,eq)_k}q^k,
\end{align}
where at least one of $a,b,c$ is of the form $q^n$, $n=1,2,\ldots.$

Furthermore, writing
$$
\frac{1}{(e/q)_k}=\frac{1}{(1-e/q)}
\left(\frac{1}{(e)_k}-\frac{1}{(e)_k}eq^{k-1}\right),
$$
we have
\begin{align*}
&\hskip -3mm
\sum_{k=-\infty}^\infty
\frac{(q/a,q/b,q/c,q/d,q/e)_k}{(a,b,c,d,e/q)_k}(abcdeq^{-3})^k \\
&=\frac{1}{1-e/q}\sum_{k=-\infty}^\infty
\frac{(q/a,q/b,q/c,q/d,q/e)_k}{(a,b,c,d,e/q)_k}(abcdeq^{-3})^k \\
&\quad -\frac{e/q}{1-e/q}\sum_{k=-\infty}^\infty
\frac{(q/a,q/b,q/c,q/d,q/e)_k}{(a,b,c,d,e/q)_k}(abcdeq^{-2})^k \\
&=\sum_{k=-\infty}^\infty
\frac{(q/a,q/b,q/c,q/d,q/e)_k}{(a,b,c,d,e/q)_k}(abcdeq^{-3})^k.
\end{align*}
The last equality holds because of the following relation
$$
\frac{(q/a,q/b,q/c,q/d,q/e)_k}{(a,b,c,d,e/q)_k}(abcdeq^{-3})^k
=\frac{(q/a,q/b,q/c,q/d,q/e)_{-k}}{(a,b,c,d,e/q)_{-k}}(abcdeq^{-2})^{-k}.
$$
Thus, by \eqref{eq:abcde1}, we have
\begin{align*}
&\hskip -3mm \sum_{k=-\infty}^\infty
\frac{(q/a,q/b,q/c,q/d,q/e)_k}{(a/q,b,c,d,e)_k}(abcdeq^{-3})^k\\
&=\sum_{k=-\infty}^\infty
\frac{(q/a,q/b,q/c,q/d,q/e)_k}{(a,b,c,d,e/q)_k}(abcdeq^{-3})^k
\\
&=\frac{(q,ab/q,bc/q,ac/q)_\infty}{(a,b,c,abc/q^2)_\infty}
\sum_{k=0}^\infty\frac{(q/a,q/b,q/c,de/q)_k}{(q,q^3/abc,d,e)_k}q^k,
\end{align*}
where at least one of $a,b,c$ is of the form $q^n$, $n=1,2,\ldots.$

Finally, setting $a=q^{n+1}$, $b=q^{l+1}$, $c=q^{m+1}$, $d=q^{u+1}$ and
$e=q^{v+1}$, we  can rewrite the last identity as follows:
\begin{align*}
&\hskip -3mm
\sum_{k=-\infty}^\infty\frac{(-1)^k q^{(5k^2-k)/2} (q)_{l+m}{(q)_{l+n}(q)_{m+n}}
(q)_{u}(q)_{v-1}(q)_{u+v}}
{(q)_{l-k}(q)_{m-k}(q)_{n-k}(q)_{u-k}(q)_{v-k}(q)_{l+k}(q)_{m+k}(q)_{n+k}
(q)_{u+k}(q)_{v+k-1}}\\
&=\sum_{k=0}^\infty\frac{q^{k^2} (q)_{l+m+n-k}(q)_{u+v+k}}
{(q)_k (q)_{l-k}(q)_{m-k}(q)_{n-k} (q)_{u+k}(q)_{v+k}}.
\end{align*}

%%%%%%%%%%%%%%%%%%%%%%%%%%%%%%%%%%%%%%%%%%%%%%%%%%%%%%%%%%%%%%%%%%%%%%%%%%
\section{Proof of \eqref{eq:abcde4}}\label{sec:abcde4}
%%%%%%%%%%%%%%%%%%%%%%%%%%%%%%%%%%%%%%%%%%%%%%%%%%%%%%%%%%%%%%%%%%%%%%%%%
We prove the equivalent form \eqref{eq:lmnrs4}.
 Writing
$(q)_{l+m+n-k}=(q)_{l+m+n-k+1}+(q)_{l+m+n-k}q^{l+m+n-k+1}$, the
left-hand side of \eqref{eq:lmnrs4} is equal to
\begin{align*}
&\hskip -3mm \sum_{k=0}^\infty\frac{q^{k^2+k}
(q)_{l+m+n-k+1}(q)_{u+v+k-1}} {(q)_k (q)_{l-k}(q)_{m-k}(q)_{n-k}
(q)_{u+k}(q)_{v+k}} +\sum_{k=0}^\infty\frac{q^{k^2+l+m+n+1}
(q)_{l+m+n-k}(q)_{u+v+k-1}} {(q)_k (q)_{l-k}(q)_{m-k}(q)_{n-k}
(q)_{u+k}(q)_{v+k}}.
\end{align*}
 Note that the first sum is equal to the left-hand side of \eqref{eq:lmnrs2}
with $u\to u-1$ and $v\to v-1$,
 while the second
sum corresponds to the left-hand side of \eqref{eq:lmnrs3}. Now
substituting the above sums by the corresponding right-hand sides
of  \eqref{eq:lmnrs2} and \eqref{eq:lmnrs3} with
 $k\to -k$ in  the second sum, and then  putting any two terms
indexed by $k$ and $-k-1$ together, the left-hand side of
\eqref{eq:lmnrs4} may be written as $\sum_{k=0}^\infty S_k$, where
\begin{align}
S_k&=\frac{(-1)^k
q^{(5k^2+3k)/2}(1-q^{2k+1})(q)_{l+m+1}(q)_{m+n+1}(q)_{l+n+1}
(q)_{u-1}(q)_{v-1}(q)_{u+v-1}}{(q)_{l-k}(q)_{m-k}(q)_{n-k}(q)_{u-k-1}(q)_{v-k-1}
(q)_{l+k+1}(q)_{m+k+1}(q)_{n+k+1}(q)_{u+k}(q)_{v+k}} \nonumber\\
&\quad + \frac{(-1)^k q^{(5k^2+k)/2+l+m+n+1}
(q)_{l+m}{(q)_{l+n}(q)_{m+n}}
(q)_{u-1}(q)_{v-1}(q)_{u+v-1}}{(q)_{l-k}(q)_{m-k}(q)_{n-k}(q)_{u-k-1}(q)_{v-k-1}
(q)_{l+k}(q)_{m+k}(q)_{n+k}(q)_{u+k}(q)_{v+k}} \nonumber\\
&\qquad\qquad\times\left(1-\frac{q^{4k+2}(1-q^{l-k})(1-q^{m-k})(1-q^{n-k})}
{(1-q^{l+k+1})(1-q^{m+k+1})(1-q^{n+k+1})}\right),  \label{eq:twosum4.2}
\end{align}
while the right-hand side of \eqref{eq:lmnrs4} may be written as
$\sum_{k=0}^\infty T_k$, where
\begin{align}
T_k&=\frac{(-1)^k q^{(5k^2+3k)/2} (q)_{l+m}(q)_{l+n}(q)_{m+n}
(q)_{u-1}(q)_{v-1}(q)_{u+v-1}}{(q)_{l-k}(q)_{m-k}(q)_{n-k}(q)_{u-k-1}(q)_{v-k-1}
(q)_{l+k}(q)_{m+k}(q)_{n+k}(q)_{u+k}(q)_{v+k}}  \nonumber\\
&\qquad\times\left(1-\frac{q^{2k+1}(1-q^{l-k})(1-q^{m-k})(1-q^{n-k})}
{(1-q^{l+k+1})(1-q^{m+k+1})(1-q^{n+k+1})}\right).  \label{eq:twosum4.3}
\end{align}
It is easy to see that $S_k=T_k$ ($k\geq 0$), which is equivalent to
\begin{align*}
&\hskip -3mm
(1-ab)(1-bc)(1-ac)(1-d^2)+d^2(1-a/d)(1-b/d)(1-c/d)(1-abcd) \\
&=(1-ad)(1-bd)(1-cd)(1-abc/d),
\end{align*}
where $a=q^{l+\frac12}$, $b=q^{m+\frac12}$, $c=q^{n+\frac12}$ and
$d=q^{k+\frac12}$. Therefore $\sum_{k\geq
0}S_k=\sum_{k\geq 0}T_k$. This completes the proof. \qed

%%%%%%%%%%%%%%%%%%%%%%%%%%%%%%%%%%%%%%%%%%%%%%%%%%%%%%%%%%%%%%%%%%%%%%%%%%%%%%%%%%%%%%%%%%%%%%%%%%
\section{Some special cases}
\label{sec:special}

In this section we give some interesting special cases and limiting
cases of Theorems~\ref{thm:lmnrs12} and \ref{thm:abcde34}. Letting
$e\to 0$ in \eqref{eq:abcde1} and \eqref{eq:abcde2},
we obtain
\begin{cor}For $|bc|<|q|$,  the following identity holds
\begin{align}
\sum_{k=-\infty}^\infty
\frac{(q/a, q/b,q/c,q/d)_k}{(a,b,c,d)_k}(-abcd)^k q^{(k^2-5k)/2}
=\frac{(q,bc/q)_\infty}{(b,c)_\infty}
\sum_{k=0}^\infty\frac{(q/b,q/c,ad/q)_k}{(q,d,a)_k}\left(\frac{bc}{q}\right)^k,
\label{eq:bcde1}
\end{align}
and for $|bc|<1$ there holds
\begin{align}
\sum_{k=-\infty}^\infty
\frac{(q/b,q/c,q/d,q/e)_k}{(bq,cq,dq,eq)_k}(-bcde)^k q^{(k^2-k)/2}
=\frac{(q,bc)_\infty}{(bq,cq)_\infty}
\sum_{k=0}^\infty\frac{(q/b,q/c,de)_k}{(q,dq,eq)_k} (bc)^k,
\label{eq:bcde2}
\end{align}
provided that at least one of $a$, $b$ and $c$ is of the form
$q^n$, $n=1,2,\ldots$.
\end{cor}

Letting $a=q^l$ and $l\to\infty$ in \eqref{eq:abcde3} and \eqref{eq:abcde4}, we
obtain
\begin{cor}For $|bc|<|q|$, there holds
\begin{align*}
\sum_{k=-\infty}^\infty
\frac{(q/b,q/c,q/d,q/e)_k}{(b,c,d/q,e/q)_k}(-bcde)^k
q^{(k^2-5k)/2} =\frac{(q,bc/q)_\infty}{(b,c)_\infty}
\sum_{n=0}^\infty\frac{(q/b,q/c,de/q^2)_k}{(q,d,e)_k}\left(\frac{bc}{q}\right)^k,
\end{align*}
and for $|bc|<1$ there holds
\begin{align*}
\sum_{k=-\infty}^\infty
\frac{(q/b,q/c,q/d,q/e)_k}{(b,c,d/q,e/q)_k}(-bcde)^k
q^{(k^2-7k)/2} =\frac{(q,bc/q)_\infty}{(b,c)_\infty}
\sum_{n=0}^\infty\frac{(q/b,q/c,de/q^2)_k}{(q,d,e)_k}(bc)^k.
\end{align*}
\end{cor}

Letting $v\to\infty$ in \eqref{eq:lmnrs1}--\eqref{eq:lmnrs4} we obtain the
following two corollaries.

\begin{cor}For $l,m,n,u\in\mathbb{N}$, there holds
\begin{align}
&\hskip -3mm
\sum_{k=0}^\infty\frac{q^{k^2} (q)_{l+m+n-k}}
{(q)_k (q)_{l-k}(q)_{m-k}(q)_{n-k} (q)_{u+k}} \nonumber\\
&=\sum_{k=-\infty}^\infty\frac{(-1)^k q^{(5k^2-k)/2} (q)_{l+m}(q)_{l+n}(q)_{m+n}(q)_u}
{(q)_{l-k}(q)_{m-k}(q)_{n-k}(q)_{u-k}(q)_{l+k}(q)_{m+k}(q)_{n+k}(q)_{u+k}},
\label{eq:lmnr1}
\end{align}
and
\begin{align}
&\hskip -3mm
\sum_{k=0}^\infty\frac{q^{k^2+k} (q)_{l+m+n-k+1}}
{(q)_k (q)_{l-k}(q)_{m-k}(q)_{n-k}(q)_{u+k}} \nonumber\\
&=\sum_{k=-\infty}^\infty\frac{(-1)^k q^{(5k^2-3k)/2} (q)_{l+m+1}(q)_{l+n+1}(q)_{m+n+1}(q)_u}
{(q)_{l-k}(q)_{m-k}(q)_{n-k}(q)_{u-k}(q)_{l+k+1}(q)_{m+k+1}(q)_{n+k+1}(q)_{u+k+1}}.
\label{eq:lmnr2}
\end{align}
\end{cor}

\begin{cor}For $l,m,n,u\in\mathbb{N}$, there holds
\begin{align}
&\hskip -3mm
\sum_{k=0}^\infty\frac{q^{k^2} (q)_{l+m+n-k}}
{(q)_k (q)_{l-k}(q)_{m-k}(q)_{n-k} (q)_{u+k}} \nonumber\\
&=\sum_{k=-\infty}^\infty\frac{(-1)^k q^{(5k^2-k)/2} (q)_{l+m}{(q)_{l+n}(q)_{m+n}}(q)_{u-1}}
{(q)_{l-k}(q)_{m-k}(q)_{n-k}(q)_{u-k}(q)_{l+k}(q)_{m+k}(q)_{n+k}(q)_{u+k-1}},
\label{eq:lmnr3}
\end{align}
and
\begin{align}
&\hskip -3mm
\sum_{k=0}^\infty\frac{q^{k^2+k} (q)_{l+m+n-k}}
{(q)_k (q)_{l-k}(q)_{m-k}(q)_{n-k}(q)_{u+k}} \nonumber\\
&=\sum_{k=-\infty}^\infty\frac{(-1)^k q^{(5k^2-3k)/2} (q)_{l+m}(q)_{m+n}(q)_{l+n}(q)_{u-1}}
{(q)_{l-k}(q)_{m-k}(q)_{n-k}(q)_{u-k}(q)_{l+k}(q)_{m+k}(q)_{n+k}(q)_{u+k-1}}.
\label{eq:lmnr4}
\end{align}
\end{cor}

Replacing $q$ by $q^{-1}$ in \eqref{eq:lmnr1}--\eqref{eq:lmnr4} and noticing that
$(q^{-1};q^{-1})_m=(-1)^mq^{-{m+1\choose 2}}(q;q)_m$, we get the following results:

\begin{cor}For $l,m,n,u\in\mathbb{N}$, there holds
\begin{align*}
&\hskip -3mm
\sum_{k=0}^{n}\frac{q^{k^2+uk} (q)_{l+m+n-k}}
{(q)_k (q)_{l-k}(q)_{m-k}(q)_{n-k} (q)_{u+k}} \nonumber\\
&=\sum_{k=-n}^{n}\frac{(-1)^k q^{(3k^2-k)/2} (q)_{l+m}(q)_{m+n}{(q)_{l+n}}(q)_u}
{(q)_{l-k}(q)_{m-k}(q)_{n-k}(q)_{u-k}(q)_{l+k}(q)_{m+k}(q)_{n+k}(q)_{u+k}},
%\label{eq:lmnr5}
\end{align*}
and
\begin{align*}
&\hskip -3mm
\sum_{k=0}^{n}\frac{q^{k^2+(u+1)k} (q)_{l+m+n-k+1}}
{(q)_k (q)_{l-k}(q)_{m-k}(q)_{n-k}(q)_{u+k+1}} \nonumber\\
&=\sum_{k=-n-1}^{n}\frac{(-1)^k q^{(3k^2+k)/2} (q)_{l+m+1}(q)_{m+n+1}(q)_{l+n+1}(q)_u}
{(q)_{l-k}(q)_{m-k}(q)_{n-k}(q)_{u-k}(q)_{l+k+1}(q)_{m+k+1}(q)_{n+k+1}(q)_{u+k+1}}.
%\label{eq:lmnr6}
\end{align*}
\end{cor}

\begin{cor}For $l,m,n,u\in\mathbb{N}$, there holds
\begin{align*}
&\hskip -3mm
\sum_{k=0}^\infty\frac{q^{k^2+uk} (q)_{l+m+n-k}}
{(q)_k (q)_{l-k}(q)_{m-k}(q)_{n-k} (q)_{u+k}} \nonumber\\
&=\sum_{k=-\infty}^\infty\frac{(-1)^k q^{(3k^2-k)/2} (q)_{l+m}{(q)_{l+n}(q)_{m+n}}(q)_{u-1}}
{(q)_{l-k}(q)_{m-k}(q)_{n-k}(q)_{u-k}(q)_{l+k}(q)_{m+k}(q)_{n+k}(q)_{u+k-1}},
%\label{eq:lmnr7}
\end{align*}
and
\begin{align*}
&\hskip -3mm
\sum_{k=0}^\infty\frac{q^{k^2+(u-1)k} (q)_{l+m+n-k}}
{(q)_k (q)_{l-k}(q)_{m-k}(q)_{n-k}(q)_{u+k}} \nonumber\\
&=\sum_{k=-\infty}^\infty\frac{(-1)^k q^{(3k^2+k)/2} (q)_{l+m}(q)_{m+n}(q)_{l+n}(q)_{u-1}}
{(q)_{l-k}(q)_{m-k}(q)_{n-k}(q)_{u-k}(q)_{l+k}(q)_{m+k}(q)_{n+k}(q)_{u+k-1}}.
%\label{eq:lmnr8}
\end{align*}
\end{cor}

Observing the symmetry of the denominators in the right-hand sides of
\eqref{eq:lmnrs1} and \eqref{eq:lmnrs2}, we obtain
\begin{cor}
For $l,m,n,u,v\in\mathbb{N}$, there holds
\begin{align}
&\hskip -3mm
\frac{1}{(q)_{l+m}(q)_{l+n}(q)_u(q)_v}
\sum_{k=0}^\infty\frac{q^{k^2} (q)_{l+m+n-k}(q)_{u+v+k}}
{(q)_k (q)_{l-k}(q)_{m-k}(q)_{n-k} (q)_{u+k}(q)_{v+k}} \nonumber\\
&=\frac{1}{(q)_{l+u}(q)_{l+v}(q)_m(q)_n}
\sum_{k=0}^\infty\frac{q^{k^2} (q)_{l+u+v-k}(q)_{m+n+k}}
{(q)_k (q)_{l-k}(q)_{u-k}(q)_{v-k} (q)_{m+k}(q)_{n+k}},         \label{eq:lmnrs5}
\end{align}
and
\begin{align}
&\hskip -3mm
\frac{1}{(q)_{l+m+1}(q)_{l+n+1}(q)_u(q)_v}
\sum_{k=0}^\infty\frac{q^{k^2+k} (q)_{l+m+n-k+1}(q)_{u+v+k+1}}
{(q)_k (q)_{l-k}(q)_{m-k}(q)_{n-k} (q)_{u+k+1}(q)_{v+k+1}} \nonumber\\
&=\frac{1}{(q)_{l+u+1}(q)_{l+v+1}(q)_m(q)_n}
\sum_{k=0}^\infty\frac{q^{k^2} (q)_{l+u+v-k+1}(q)_{m+n+k+1}}
{(q)_k (q)_{l-k}(q)_{u-k}(q)_{v-k} (q)_{m+k+1}(q)_{n+k+1}}.    \label{eq:lmnrs6}
\end{align}
\end{cor}

Letting $l,u,v\to\infty$ in \eqref{eq:lmnrs5} and \eqref{eq:lmnrs6}, we get
\begin{align}
\frac{1}{(q)_\infty}\sum_{k=0}^\infty\frac{q^{k^2}}{(q)_k (q)_{n-k}(q)_{m-k}}
&=\frac{1}{(q)_m(q)_n}
\sum_{k=0}^{\infty}\frac{q^{k^2}(q)_{m+n+k}}{(q)_k (q)_{m+k}(q)_{n+k}}, \label{eq:eulermn1} \\
\frac{1}{(q)_\infty}\sum_{k=0}^\infty\frac{q^{k^2+k}}{(q)_k (q)_{n-k}(q)_{m-k}}
&=\frac{1}{(q)_m(q)_n}
\sum_{k=0}^{\infty}\frac{q^{k^2+k}(q)_{m+n+k+1}}{(q)_k (q)_{m+k+1}(q)_{n+k+1}},
\label{eq:eulermn2}
\end{align}
Furthermore, letting $m\to\infty$ in \eqref{eq:eulermn1} and \eqref{eq:eulermn2},
we obtain the following remarkable identities:
\begin{align}
\frac{1}{(q)_\infty}\sum_{k=0}^\infty\frac{q^{k^2}}{(q)_k (q)_{n-k}}
&=\frac{1}{(q)_n} \sum_{k=0}^{\infty}\frac{q^{k^2}}{(q)_k (q)_{n+k}},
   \label{eq:eulern1} \\
\frac{1}{(q)_\infty}\sum_{k=0}^\infty\frac{q^{k^2+k}}{(q)_k (q)_{n-k}}
&=\frac{1}{(q)_n}\sum_{k=0}^{\infty}\frac{q^{k^2+k}}{(q)_k (q)_{n+k+1}}.  \nonumber
%\label{eq:eulern2}
\end{align}

In what follows we assume that $\frac{1}{n!}=0$ if $n<0$. Letting
$q$ tend to $1$ in \eqref{eq:lmnrs1}, we obtain
\begin{cor}For $l,m,n,u,v\in\mathbb{N}$, there holds
\begin{align*}
&\hskip -3mm
\sum_{k=-\infty}^\infty
(-1)^k{l+m\choose l+k}{m+n\choose m+k}{n+l\choose n+k}{u+v\choose u+k}{u+v\choose v+k} \\
&={u+v\choose u}\sum_{k=0}^{\infty}
\frac{(l+m+n-k)!(u+v+k)!}{k!(l-k)!(m-k)!(n-k)!(u+k)!(v+k)!}.
\end{align*}
\end{cor}
In particular, when $l=m=n=u=v$, we get
\begin{align}
\sum_{k=-n}^n (-1)^k{2n\choose n+k}^5
={2n\choose n}\sum_{k=0}^n {3n-k\choose n-k}{2n+k\choose k}{2n\choose n+k}^2.
\label{eq:bino5}
\end{align}
Replacing $k$ by $n-k$ in the right-hand side of \eqref{eq:bino5} yields
$$
\sum_{k=-n}^n (-1)^k{2n\choose n+k}^5
={2n\choose n}\sum_{k=0}^n {3n-k\choose n-k}{2n+k\choose k}{2n\choose k}^2.
$$

Similarly, letting $q$ tend to $1$ in \eqref{eq:lmnr1}, or first letting
$l$ tend to infinity in \eqref{eq:lmnrs1} and then letting $q$ tend to $1$, we obtain
the following formulae:
\begin{cor}For $l,m,n,u,v\in\mathbb{N}$, there holds
\begin{align*}
&\hskip -3mm
\sum_{k=-\infty}^\infty
(-1)^k{l+m\choose l+k}{m+n\choose m+k}{n+l\choose n+k}{2u\choose u+k} \\
&=\frac{(2u)!}{u!}\sum_{k=0}^{\infty}
\frac{(l+m+n-k)!}{k!(l-k)!(m-k)!(n-k)!(u+k)!}, %\label{eq:bino4-1}
\end{align*}
and
\begin{align*}
&\hskip -3mm
\sum_{k=-\infty}^\infty
(-1)^k{m+n\choose m+k}{m+n\choose n+k}{u+v\choose u+k}{u+v\choose v+k} \\
&={u+v\choose u}\sum_{k=0}^{\infty}
\frac{(m+n)!(u+v+k)!}{k!(m-k)!(n-k)!(u+k)!(v+k)!}.    %\label{eq:bino4-2}
\end{align*}
\end{cor}
In particular, when $l=m=n=u=v$, we get
\begin{align}
\sum_{k=-n}^n (-1)^k{2n\choose n+k}^4
&={2n\choose n}\sum_{k=0}^n {3n-k\choose n-k}{2n\choose n+k}{n\choose k} \nonumber\\
&={2n\choose n}\sum_{k=0}^n {2n+k\choose k}{2n\choose n+k}^2. \label{eq:bino4}
\end{align}
From \eqref{eq:bino5} and \eqref{eq:bino4} one sees that
$\sum_{k=-n}^n (-1)^k{2n\choose n+k}^m$ is divisible by
${2n\choose n}$ for $m=4,5$. In a separate paper \cite{GZfactor}, we have generalized
this result by showing that the alternating sum
$$
\sum_{k=-n}^n (-1)^k\prod_{i=1}^m{n_i+n_{i+1}\choose n_i+k}\qquad (n_{m+1}=n_1)
$$
is a nonnegative integer divisible by all the binomial coefficients
${n_j+n_{j+1}\choose n_j}$ for $1\leq j\leq m$.

\section*{Acknowledgments} This work was partially done during the first
author's visit to Institut Camille Jordan, Universit\'e
Claude Bernard (Lyon I), and supported by a French postdoctoral fellowship.
The first author was also supported by Shanghai Leading Academic Discipline Project B407.
The authors thank the two anonymous referees for useful comments on a previous version of this paper.

\renewcommand{\baselinestretch}{1}

\end{document}